\documentclass{article}
\usepackage{setspace}
\usepackage{graphicx} % Required for inserting images
\usepackage[top=1in, bottom=1.25in, left=1.25in, right=1.25in]{geometry}
\usepackage{mathtools}
\usepackage[toc,page]{appendix}
\usepackage{enumerate}
\usepackage{amsmath}
\usepackage{empheq}
\usepackage{amssymb}
\usepackage{amsthm}
\usepackage{listings} 
\usepackage{xcolor} 
\usepackage{graphicx} 

\usepackage{caption}
\usepackage{subcaption}
\usepackage{enumerate}
\usepackage{color}
\usepackage{float}
\usepackage{csvsimple}
\usepackage{comment}
\usepackage[hidelinks]{hyperref}
\usepackage{cleveref}
\usepackage{framed}
\usepackage{mwe}
\usepackage{placeins} 
\theoremstyle{remark}
\newtheorem{remark}{Remark}

\begin{document}
\title{An anisotropic traffic flow model with look-ahead effect for mixed autonomy traffic}
\author{Shouwei Hui \thanks{Shouwei Hui is with the Department
of Mathematics, Univercity of California Davis, Davis,
CA, 95616 USA e-mail: (huihui@ucdavis.edu).}, Michael Zhang \thanks{Michael Zhang is with the Department
of Civil and Environmental Engineering, Univercity of California Davis, Davis,
CA, 95616 USA e-mail: (hmzhang@ucdavis.edu)} \thanks{Corresponding author. Accepted by TRB Annual Meeting 2025}}
\date{}
% <-this % stops a space

\maketitle
\section*{Abstract}

In this paper we extend the Aw-Rascle-Zhang (ARZ) non-equilibrium traffic flow model to take into account the look-ahead capability of connected and autonomous vehicles (CAVs), and the mixed flow dynamics of human driven and autonomous vehicles. The look-ahead effect of CAVs is captured by a non-local averaged density within a certain distance (the look-ahead distance). We show, using wave perturbation analysis, that increased look-ahead distance loosens the stability criteria. Our numerical experiments, however, showed that a longer look-ahead distance does not necessarily lead to faster convergence to equilibrium states. We also examined the impact of spatial distributions and market penetrations of CAVs and showed that increased market penetration helps stabilizing mixed traffic while the spatial distribution of CAVs have less effect on stability. The results revealed the potential of using CAVs to stabilize traffic, and may provide qualitative insights on speed control in the mixed autonomy environment.

% \hfill\break%
% \noindent\textit{Keywords}: Connected and autonomous vehicles, look-ahead, wave-perturbation analysis, numerical experiments
% \newpage

\section{Introduction and Review of Related Work} \label{sec1}
% Literature review on macroscopic models, mixed flow control, etc 6 papers on PDE, 3 to 4 papers on multiclass wong wong mohan ramadurai goatin hamori
% LWR instantly move to equilibrium state which is unrealistic.
% ARZ relaxation is more reasonable, overcome previous flaws, but not exactly mimicing HDVs, GSOM
% with macroscopic models: Non-local models Chiarello, Hamori, finite time shock formation in Hamori
% Multiclass: Wong, Ngoduy*2, Logghe, Zhang, Qian, Mohan, Huang moving bottleneck as controller; (PDE-ODE Piccoli, Goatin, Monache)
% CAVs: With development of technology, CAVs as controller: control models on CAVs: , Yu's boundary control, Huang's mixed flow CF models mention Zhu, Wang, Zhao and find some more in
% Field experiments: 3 in total

\subsection{Hydrodynamic traffic flow models}
%LWR +2
%PW +2
%ARZ +2
%GSOM +1 Lebacque
%Nonlocal and shock formation Chiarello, Hamori, Ramadan, 2 extra +5
Hydrodynamic traffic flow models, often given in the form of partial differential equations (PDEs) have been widely studied in the traffic science literature. They have wide applications and are often used in traffic simulation, state estimation and control design. The most classic of them is the Lighthill-Whitham-Richards (LWR) model \cite{lighthill1955kinematic, richards1956shock}, which has the form
\begin{equation} 
       \rho_t+(\rho V(\rho))_x=0,
       \label{lwr}
   \end{equation}
where $\rho$ is the density of traffic at location $x$ and time $t$, and $V(\rho)$ is the equilibrium velocity as a function of density. The LWR model is essentially a scalar conservation law endowed with a equation of state that captures average driver behavior under stationary (or equilibrium) conditions. The LWR model is capable of modeling transitions from one stationary state to another, in the form of shock or acceleration waves. However, it lacks the ability to model some notable traffic flow phenomena, such as stop-and-go waves and traffic hysteresis. Various models, collectively known as higher-order traffic flow models, have been proposed to overcome LWR model's deficiencies. For example, analogous to shallow channel water flows, Payne and Whitham \cite{payne1971model, whitham2011linear} respectively introduced a momentum equation to
capture speed evolution away from equilibrium, and proposed the first higher-order model:

\begin{align}
    \begin{cases}
           \rho_t+(\rho v)_x=0,\\
       v_t+vv_x=\frac{V(\rho)-v}{\tau}-\frac{c_0^2}{\rho}\rho_x,
    \end{cases}\label{pw}
    \end{align}
   where $\tau$ is a relaxation time constant and $c_0<0$ is the "traffic sound speed". But this model has two main drawbacks: it can produce negative travel speed ('wrong way travel') and traffic information can travel faster than vehicles, which violates the anisotropic property of traffic flow---that is, vehicles cannot push other vehicles from behind to speed them up. To solve these problems, there are two models independently introduced in \cite{aw2000resurrection,zhang2002non}, and the inhomogeneous Aw-Rascle-Zhang (ARZ) model has the form:
   \begin{align}
   \begin{cases}
   \rho_t+(\rho v)_x=0,\\
       (v+h(\rho))_t+v(v+h(\rho))_x=\frac{V(\rho)-v}{\tau},
   \end{cases}\label{arz}
   \end{align}
   where the constant $c_0$ in PW model is substituted by the convective derivative ($\partial_t+v\partial_x$) of the pressure function $h(\rho)$ accounting for drivers’ anticipation of downstream density changes.

   The ARZ model has been widely used and studied since it was first proposed. Theoretical and numerical solutions of the ARZ model have been studied in  e.g. \cite{lebacque2007aw, mammar2009riemann, mohammadian2018improved}. Others have extended the ARZ model: for example,  Lebacque et al \cite{lebacque2007generic} generalised the ARZ model to the generic second order models (GSOM) where the pressure term is generalised to a non-linear velocity term. The GSOM model have then been used for data fitting \cite{fan2013comparative} and extended with non-local densities \cite{chiarello2020micro, hamori2024aw}. 

\subsection{Multi-class hydrodynamic traffic flow models}
% First order: wong, Zhang, Logghe, Ngoduy, Qian +5
% Second order: Ngoduy, Mohan, Huang +3
Real-world traffic has vehicles of different types and performances, which can be categorized into vehicle classes. Each class of vehicles may interact with others in different ways and this can be captured by extending the aforementioned models to multi-class hydrodynamic traffic flow models. Starting with an extension of the LWR model, Wong and Wong \cite{wong2002multi} proposed a multi-class LWR model with heterogeneous drivers characterized by their choice of free-flow speeds. In particular, they gave an isotropic case where the speed of each class is a function of the total density. In a separate work, Zhang and Jin \cite{zhang2002kinematic} proposed a multi-class LWR model considering critical density such that when traffic concentration reached a critical value, all the class of vehicles are mixed together and move as a group, and below the critical density the model is similar to Wong and Wong's model. Ngoduy and Liu \cite{ngoduy2007multiclass} proposed a generalized multi-class first-order simulation model based on an approximate Riemann solver, which is able to explain certain non-linear traffic phenomena on freeways. Logghe and Immers \cite{logghe2008multi} proposed a new model where vehicle classes interact in a non-cooperative way, where slow vehicles act as moving bottlenecks for fast vehicles while fast vehicles have no influence on slow vehicles. Such relations have been previous presented in \cite{newell1998moving}. Qian et al \cite{qian2017modeling} developed a macroscopic heterogeneous traffic flow model with pragmatic cross-class interaction rules.

There are also studies that proposed non-equilibrium hydrodynamic models for mixed traffic flow, e.g. \cite{ngoduy2009continuum, ngoduy2013analytical, ngoduy2013instability, mohan2017heterogeneous, huang2020scalable}. Specifically, Ngoduy et al \cite{ngoduy2009continuum} proposed a multi-class gas-kinetic model where one class of vehicles are able to receive a warning massage when there is downstream congestion and further extended it in \cite{ngoduy2013analytical, ngoduy2013instability} to include cooperative adaptive cruise control (CACC). Mohan and Ramadurai \cite{mohan2017heterogeneous} extends the ARZ model to a multi-class model using area occupancy (AO) which can capture the unique phenomena in lane-free traffic. Huang et al \cite{huang2020scalable} proposed a multi-class model where human driven vehicles (HDVs) are modeled by the ARZ model and CAVs are modeled by a mean-field game. They also performed linear stability analysis for the mean-field game model. 

\subsection{CAVs as agents for traffic stabilization}
% PDE control +4
% CF control +3
% Experiment +3

Traffic flow of HDVs can be unstable even without an external disturbance. For example, in \cite{sugiyama2008traffic}, a field experiment on a ring road with human driven vehicles showed that stop-and-go waves can arise without the presence of any bottlenecks when there are sufficient number of vehicles on the road. A recent field experiment, on the other hand, showed that such stop-and-go waves can be eliminated with a single AV (Autonomous Vehicle) as a control agent to pace HDV traffic for the vehicles involved \cite{stern2018dissipation}. Such improvements were also found in a larger field experiment of over $100$ CAVs \cite{lee2024traffic}. This stabilization effect of an AV or CAVs as a control agent has also been widely studied through traffic simulation using microscopic car-following models, e.g. \cite{cui2017stabilizing, wang2021leading, zhao2023safety}. In these studies, it is shown that a single AV can stabilize multiple HDVs on a single-lane road by using its sensing capabilities and feedback control to adjust its speed.
\subsection{The main contributions of this paper}
In this paper, we enhance the understanding on the look-ahead effect of CAVs in traffic flow modelling by extending the ARZ model with a non-local density parameter, which simulates the forward-looking capabilities of CAVs. This modification allows for a more realistic representation of how autonomous technologies might influence traffic flow dynamics.

We undertake a comprehensive theoretical stability analysis using wave perturbation methods and demonstrate that the extended model for CAVs can achieve greater stability over longer look-ahead distances, offering a theoretical foundation for integrating CAVs into traffic systems.

Additionally, we further extend our model to a multi-class framework, accommodating both HDVs and CAVs. This extension is crucial to evaluate the stabilization effect of CAVs in various traffic conditions with presence of HDVs. Through extensive simulations referencing the studies above, we evaluate how different configurations of look-ahead distances and vehicle distributions impact traffic flow stability. 

The findings of this study contribute to the ongoing discussions on traffic management in the mixed autonomy environments. One of them suggests that moderate look-ahead distances might provide optimal stability conditions. Another notable finding is that with a relatively low penetration rate of CAVs, the mixed flow can be effectively stabilized, which is consistent with previous studies. Furthermore, evenly distributed CAVs achieve marginally better stabilization results compared to segregated distributions. 
% Organizaiton

The remainder of the chapter is organized as follows: Section \ref{sec32} introduces the modified ARZ model and interprets it as a model for CAVs. Section \ref{sec33} gives a stability analysis of the model using wave perturbation. Section \ref{sec34} formulates a multi-class extension of the modified ARZ model for mixed CAV-HDV traffic and in Section \ref{sec35} parameters of both models are analysed via numerical experiments to test the stability of CAVs under different conditions. Lastly in Section \ref{sec36} conclusion is drawn and directions of future research are proposed.

\section{An extended ARZ model with look-ahead effect} \label{sec32}
% Well-definedness for a similar model already been proved by Hamori 2024. Need to check his previous paper and see if need to cite

\begin{figure}[ht]
    \centering
    \includegraphics[width=0.8\linewidth]{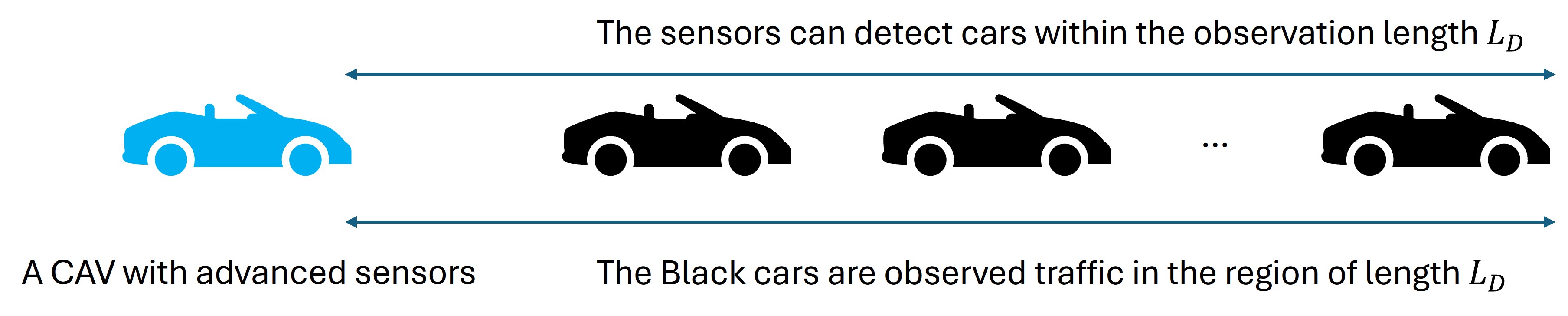}
    \caption{A CAV's front observation of the traffic density of a certain distance in front.}
    \label{CAVobserve}
\end{figure}

We first extend the ARZ model to take into account the look-ahead capability of CAVs without explicitly modeling CAVs and HDVs as distinct classes. Here we assume that the CAVs are all equipped with range sensors and vehicle to vehicle communication to enable them to observe the density of a certain distance ahead, say $L_D$. A visualised demonstration of this is given in \textbf{Figure \ref{CAVobserve}}. 

Instead of responding to the motion of the immediate vehicle in front, a CAV can take advantage of this look-ahead capability and adopt a speed that is based on the average traffic condition within this look-ahead distance, therefore reducing over- or under- reaction and smoothing its trajectory. This will in turn lead to greater stability of traffic. Following this argument, we modify the Aw-Rascle-Zhang model with a new relaxation term that takes into account this look-ahead effect on traffic flow as follows:
\begin{equation}
\begin{cases}
    \begin{aligned}
        &\rho_t+(\rho v)_x=0,\\
        &(v+h(\rho))_t+v(v+h(\rho))_x=\frac{v-V(\rho^*)}{\tau},\\
        &\rho^*(L_D)=\frac{\int_x^{x+L_D}\rho(t,\xi)d\xi}{L_D},  
    \end{aligned}
    \end{cases}
    \label{modARZ}
\end{equation}
where the relaxation of $v$ is toward an equilibrium speed $V(\rho^*)$ with $\rho^*$ as the average traffic density in the observation region $[x,x+L_D]$. The observed average density is calculated by integration. Moreover, if $L_D$ goes to $0$ the model reduces to the original ARZ model. We can also generalize the average density with a weight function for biased observation:
\begin{remark}
    A more general weighted average density $\rho^*(L_D,w)$ with weight function $w(x)$ can be defined as
\begin{equation}
    \rho^*(L_D,w)=\frac{\int_x^{x+L_D}w(\xi)\rho(t,\xi)d\xi}{L_D},  
\end{equation}
where $w(x)$ satisfies $\int_x^{x+L_D}w(\xi)d\xi=1/L_D$. With the weighted density, CAVs are considering vehicles in front with different sensitivities, similar to the microscopic multi-following model in \cite{lenz1999multi}.

With the look-ahead (weighted) average density we are primarily focusing on the CACC logic in CAVs. There are many other complex dynamics and controls which can be implemented into the model \eqref{modARZ}.

\end{remark}
Additionally for periodic boundary conditions (i.e. traffic on a ring road), partial observation (look-ahead) is equivalent to full observation (look-ahead of the entire road) when $L_D=L$, the length of the ring road. 

For readers who are interested in the theoretical analysis such as solution existence, this model can be implicitly written as the non-local traffic model in \cite{hamori2024aw} that is proven well defined under certain constraints. %With similar thoughts, the pressure function of CAVs can also be modified with a non-local density term, which means CAVs can adapt to a wider range of surrounding traffic compared to HDVs. We will consider such extensions in future work.
% In the numerical experiments we just adopt the same pressure function as in Ramadan2019. 

% Interaction follows Huang 2021, form follows Mohan 2017 but not using laneless model
\section{Stability analysis of the extended ARZ model} \label{sec33}
In this section we will follow the classic wave perturbation analysis approach \cite{kerner1993cluster}\cite{flynn2009self}\cite{ramadan2021structural} to analyze the stability of the extended ARZ model \eqref{modARZ}.

For a given initial state $(\rho_0,v_0)$, the steady state solution of the ARZ model is $(\rho,v)=\left(\rho_0, V(\rho_0)\right)$ for some $0< \rho_0 < \rho_{j}$ where $\rho_{j}$ is the jam density. Now assume that the initial condition is perturbed by a small periodic disturbance:
\begin{equation}
    \rho=\rho_0+\tilde{\rho};\;v=V(\rho_0)+\tilde{v},
\end{equation}
where
\begin{equation}
    \tilde{\rho}=Re^{ikx+\sigma t};\;\tilde{v}=Ve^{ikx+\sigma t}
\end{equation}
with $R,V$ has infinitesimal constant scales, and $k,\sigma$ are constants for the perturbation's frequency and amplitude, respectively. 

By neglecting second or higher order terms of $R$ and $V$ we can derive a linear system
\begin{equation}
    \begin{bmatrix}
        \sigma +ik\psi & ik\rho_0 \\
        \sigma\phi+ik\psi\phi-\frac{\zeta}{\tau} & \sigma+ik\psi+\frac{1}{\tau}
    \end{bmatrix}
    \begin{bmatrix}
        R \\
        U
    \end{bmatrix}
    =
    \begin{bmatrix}
        0 \\
        0
    \end{bmatrix}
\end{equation}
where $\psi=V(\rho_0)>0$, $\phi=h'(\rho)>0$, $\zeta=(e^{ikL_D}-1)V'(\rho)/(ikL_D)$. 

Follow the calculation process in \cite{ramadan2021structural}, we can deduce that traffic is stable when
\begin{equation}
    h'(\rho)+\frac{|\sin(kL_D)|}{kL_D}V'(\rho)>0
\end{equation}
Since $|\sin x|/x$ is a decreasing function of $x$, this equation implies that with certain level of oscillation frequency, stability criteria does not depend on $\tau$ and the range of stability can be increased with $L_D$.
\begin{remark}
    In particular if $\rho(x+L)=\rho(x)$ for all $x\in\mathbb{R}$, then if $L_D=L$ we have full observation of the road and the model will be always stable. In this case $\rho^*=\rho_0$ which implies that all the vehicles are relaxing toward equilibrium speed.
\end{remark}

\section{A multiclass extension of the ARZ model with look-ahead effect} \label{sec34}
In this section, we propose a model for mixed autonomy traffic where HDVs and CAVs are modeled as distinctive classes. Similar to \cite{huang2020scalable}, in our model the HDVs are reacting to total density of traffic at its position. If we let $\rho^h$ denote density of HDVs and $\rho^c$ denote the density of CAVs, then the model has the form

\begin{subequations}
\label{multiARZ}
\begin{empheq}[left=\empheqlbrace]{align}
        &\rho^h_t+(\rho^h v^h)_x=0, \\
        &\left(v^h+h(\rho^s)\right)_t+v^h\left(v^h+h(\rho^s)\right)_x=\frac{v^h-V(\rho^s)}{\tau}, \\
        &\rho^c_t+(\rho^c v^c)_x=0, \\
        &\left(v^c+h(\rho^s)\right)_t+v^c\left(v^c+h(\rho^s)\right)_x=\frac{v^c-V(\rho^*)}{\tau}, \\
        & \rho^s=\rho^h+\rho^c,
\end{empheq}    
\end{subequations}
where $v^h$ is the speed of HDVs and $v^c$ is the speed of CAVs. To highlight the look-ahead effect, for the CAVs we assume that they have the same pressure function and relaxation constant as HDVs. With similar reasons we assume that CAVs and HDVs follow the same FD. For such mixed flow the stability can depend on the ratio and distribution of vehicles, and the control method of CAVs, which means that it is hard to obtain the stability condition analytically for the traffic flow model given in \eqref{multiARZ}. In this paper, we resort to numerical solutions of \eqref{multiARZ} to explore the stability properties of this multi-class non-equilibrium model, which will be presented in the next section.
\begin{remark}
    Practically, in mixed autonomy CAVs might be capable to observe both density and speed of surrounding HDVs to change their speed accordingly, which means the pressure term and relaxation term can be defined with consideration of the density of HDVs. We will consider such extensions in future work.
\end{remark}
% since the stable solution for such multi-class traffic model may not be unique
\section{Numerical solutions} \label{sec35}
% Different waves with different traffic composition
In order to obtain numerical solutions for \eqref{modARZ} and \eqref{multiARZ}, we adapted a forward scheme with an approximate Riemann solver in \cite{ramadan2021structural} that has low computation cost and preserves properties of finite volume methods. To calculate the average density, we use a Riemann sum to give an estimation of the integration term. Given $\Delta t$, $\Delta x$ as time and space step size, $q=\rho(v+h(\rho))$ as a conserved flux variable, $i,n$ as space step variable and time step variable, and suppose that $L_D/\Delta t$ is a non-negative integer, then the update rule for approximate solutions of \eqref{modARZ} can be written as
\begin{subequations}
\begin{empheq}[left=\empheqlbrace]{align}
    &\rho_i^{n+1} =\rho_i^n - \frac{\Delta t}{\Delta x}\left((F_\rho)_{i+\frac{1}{2}}^{n}-(F_\rho)_{i-\frac{1}{2}}^{n}\right) \\
    &\rho_{i}^*=\frac{\sum_{j=1}^{L_D/\Delta t}\rho_j^{n+1}}{L_D/\Delta t} \\
    &q_i^{n+1} =q_i^n - \frac{\Delta t}{\Delta x}\left((F_q)_{i+\frac{1}{2}}^{n}-(F_q)_{i-\frac{1}{2}}^{n}\right)-\frac{\Delta t}{\tau}\left(V(\rho_i^*)+h(\rho_i^{n+1})\right)
\end{empheq}
\label{numscheme}
\end{subequations}
where the update of the approximated flow $q_i^{n+1}$ is calculated after the update of the approximated density $\rho_i^{n+1}$, and the update of the relaxation term adopts an implicit scheme to improve numerical stability. The average density $\rho_i^*$ is calculated by Riemann sum and the numerical fluxes $(F_\rho)_{i+\frac{1}{2}}^{n}, (F_q)_{i+\frac{1}{2}}^{n}$ are calculated by the Harten, Lax and van Leer (HLL) approximate Riemann solver \cite{harten1997high}. The update rule for approximate solutions of \eqref{multiARZ} can be similarly written by separately updating solutions of HDVs and CAVs using \eqref{numscheme}. 

For the model parameter values, we assumed that vehicles are on a ring road with length $L=1000$ meters, and set $\Delta t=0.05$s, $\Delta x=5$m, $\tau=3$s, $h(\rho)=8*((\rho-10)/(140-\rho))^{1/2}$m/s, similar to \cite{huang2020scalable}. The equilibrium speed model is a smooth function that combined the Greenshields model \cite{greenshields1935study} and the Triangular FD model \cite{daganzo2008analytical}: 
\begin{equation}
    V(\rho)=\begin{cases}
        \begin{aligned}
            &v_f, & \text{ if } \rho\leq\rho_f; \\
            &v_f\left(1-\frac{\rho-\rho_f}{\rho_j-\rho_f}\right), & \text{ if } \rho_f\leq\rho\leq\rho_j; \\
            &0,  & \text{ if } \rho\geq\rho_j,
        \end{aligned}
    \end{cases}
\end{equation}
where $\rho_f=10$veh/km is the free flow density, $\rho_j=140$veh/km is the jam density and $v_f=20$m/s is the free flow speed.
The initial density is a sinusoidal wave perturbation of equilibrium state similar to \cite{huang2020scalable} as well:
\begin{equation}
    \rho_0(x)=0.4*\rho_j+0.1*\rho_j*\sin(2\pi x/L),
\end{equation}
where for mixed flow we substitute $\rho$ by $\rho^s$. The initial velocity is then set as $v_0(x)=V(\rho_0(x))$. In the following subsections we will use two cases to evaluate the asymptotic stability of both models under different $L_D$ and vehicle mixes.

\subsection{Investigation of the look-ahead effect}
\begin{figure}[ht]
    \centering
        \centering
    \begin{subfigure}[b]{0.48\textwidth}
        \centering
        \includegraphics[height=0.2\textheight]{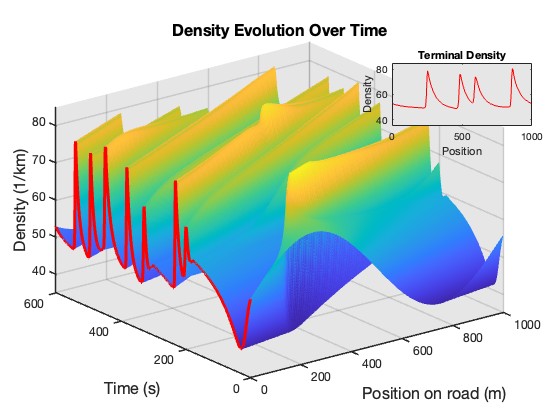}
        \caption{Density evolution}
        \label{11-1}
    \end{subfigure}
    \begin{subfigure}[b]{0.48\textwidth}
        \centering
        \includegraphics[height=0.2\textheight]{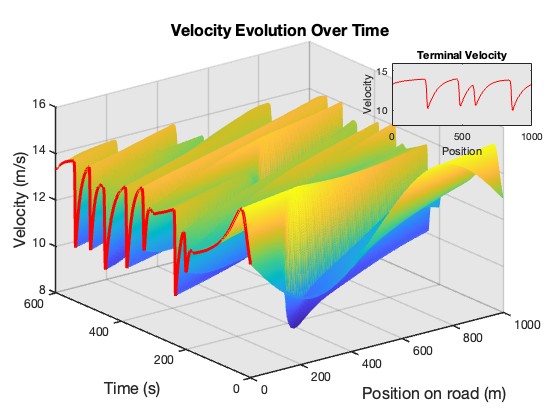}
        \caption{Velocity evolution}
        \label{11-2}
    \end{subfigure}
    \caption{Density and velocity evolution of the ARZ model ($L_D\to 0^+$), where the flow is not stable.}
    \label{SimuARZ}
\end{figure}

\begin{figure}[ht]
    \centering
        \centering
    \begin{subfigure}[b]{0.48\textwidth}
        \centering
        \includegraphics[height=0.2\textheight]{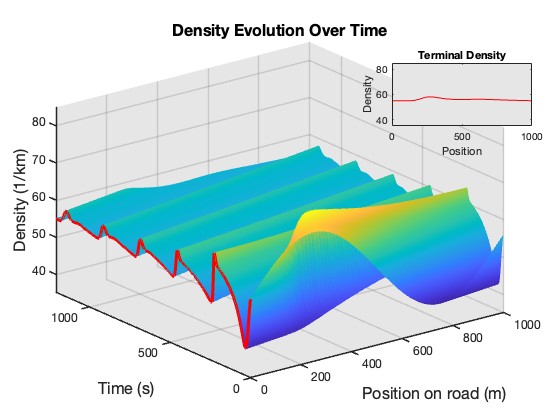}
        \caption{Density evolution}
        \label{12-1}
    \end{subfigure}
    \begin{subfigure}[b]{0.48\textwidth}
        \centering
        \includegraphics[height=0.2\textheight]{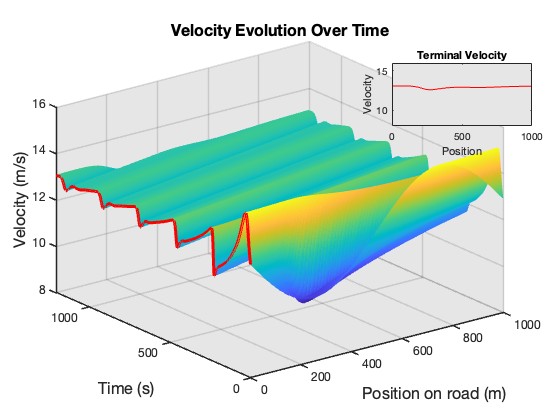}
        \caption{Velocity evolution}
        \label{12-2}
    \end{subfigure}
    \caption{Density and velocity evolution of the modified model with $L_D=15$m.}
    \label{SimuLd15}
\end{figure}

\begin{figure}[ht]
    \centering
        \centering
    \begin{subfigure}[b]{0.48\textwidth}
        \centering
        \includegraphics[height=0.2\textheight]{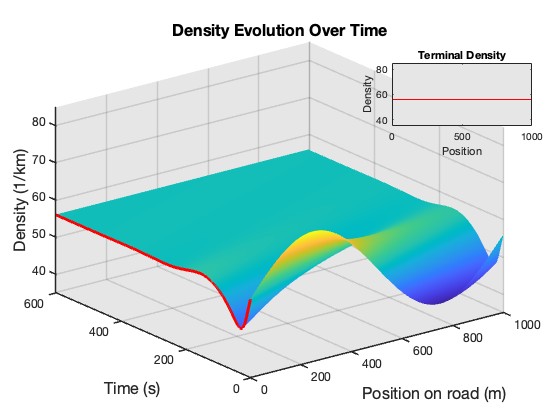}
        \caption{Density evolution}
        \label{13-1}
    \end{subfigure}
    \begin{subfigure}[b]{0.48\textwidth}
        \centering
        \includegraphics[height=0.2\textheight]{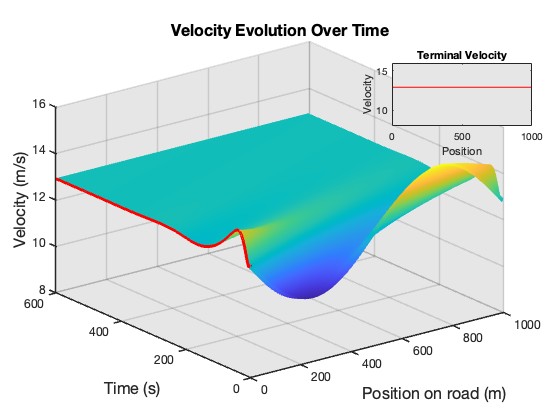}
        \caption{Velocity evolution}
        \label{13-2}
    \end{subfigure}
    \caption{Density and velocity evolution of the modified model with $L_D=100$m.}
    \label{SimuLd100}
\end{figure}

\begin{figure}[ht]
    \centering
        \centering
    \begin{subfigure}[b]{0.48\textwidth}
        \centering
        \includegraphics[height=0.2\textheight]{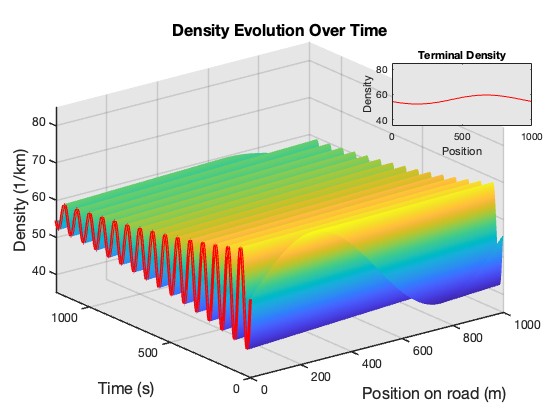}
        \caption{Density evolution}
        \label{14-1}
    \end{subfigure}
    \begin{subfigure}[b]{0.48\textwidth}
        \centering
        \includegraphics[height=0.2\textheight]{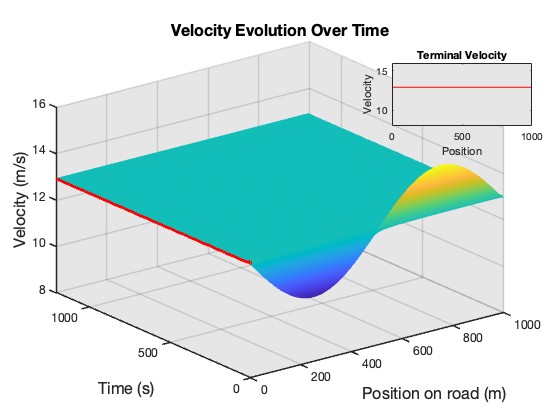}
        \caption{Velocity evolution}
        \label{14-2}
    \end{subfigure}
    \caption{Density and velocity evolution of the modified model with $L_D=1000$m.}
    \label{SimuLdmax}
\end{figure}

In this scenario we evaluate the look-ahead distance $L_D$ on the convergence of the extended ARZ model. We will consider $L_D=15,100,1000$m and compare the model results with those from the ARZ model ($L_D\to 0^+$). For $L_D=0^+,100$m we set the time duration as $T=600$s and the others as $T=1200$s. The numerical results of density and velocity evolution are shown in \textbf{Figures \ref{SimuARZ}-\ref{SimuLdmax}}.

% With too long of observation the convergence can be slower, since the speed of vehicles get too fast to equilibrium speed and not allowing the density to relax. A proper distance of observation help the most under the simulation setting
From the numerical results, we can observe that in all cases look-ahead helps stabilizing traffic, as the only unstable case is when $L_D\to 0^+$. However, longer look-ahead distance is not equivalent to faster convergence to equilibrium state. With full observation, i.e. $L_D=1000$m, or a shorter partial observation ( $L_D=15$m ) the convergence speed is much slower than with $L_D=100$m. In the case of $L_D=15$m, the look-ahead effect is not significant since this is no better than follow one-vehicle ahead. With the much longer $L_D$, the redundant information from far away is also built into drivers' response, and hence can be detrimental rather than beneficial to traffic stability when traffic conditions vary significantly over space. There seems to be a theoretical optimal look-ahead distance for achieving greater convergence and traffic stability, which may depend on parameter settings and even initial and boundary traffic conditions. We will explore this problem in our future work.   
\subsection{Investigation of stability in mixed autonomy traffic}
% mixed CAVs: two observation abilities and HDVs with different ratio, assuming that they are evenly mixed
In this scenario we investigate the potential of using CAVs to smooth and stabilize mixed traffic flow, considering two different spatial distributions of CAVs in the traffic mix.
\subsubsection{Even distribution}
\begin{figure}[ht]
    \centering
        \centering
    \begin{subfigure}[b]{0.48\textwidth}
        \centering
        \includegraphics[height=0.2\textheight]{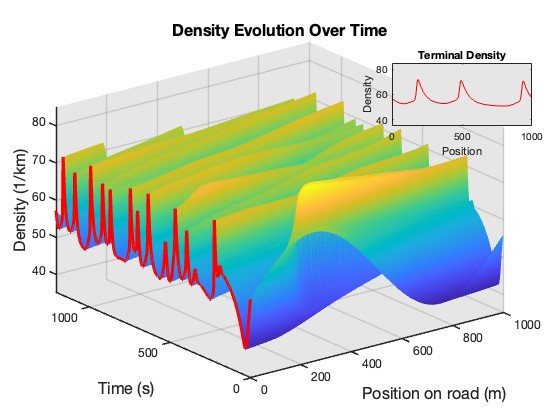}
        \caption{Density evolution}
        \label{211-1}
    \end{subfigure}
    \begin{subfigure}[b]{0.48\textwidth}
        \centering
        \includegraphics[height=0.2\textheight]{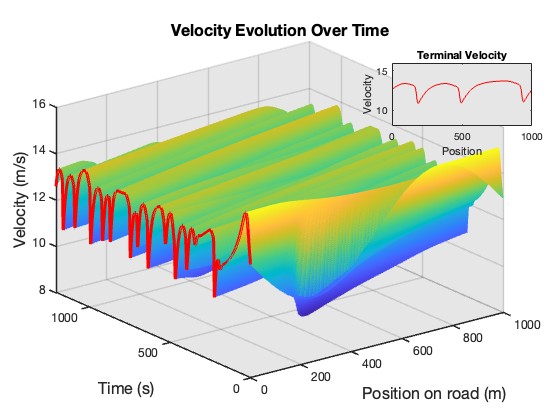}
        \caption{Velocity evolution}
        \label{211-2}
    \end{subfigure}
    \caption{Density and velocity evolution of the mixed flow model with $10\%$ of CAVs  evenly distributed.}
    \label{Simu_even10}
\end{figure}

\begin{figure}[ht]
    \centering
        \centering
    \begin{subfigure}[b]{0.48\textwidth}
        \centering
        \includegraphics[height=0.2\textheight]{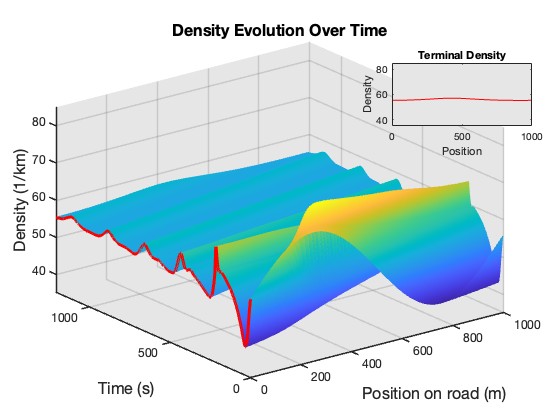}
        \caption{Density evolution}
        \label{212-1}
    \end{subfigure}
    \begin{subfigure}[b]{0.48\textwidth}
        \centering
        \includegraphics[height=0.2\textheight]{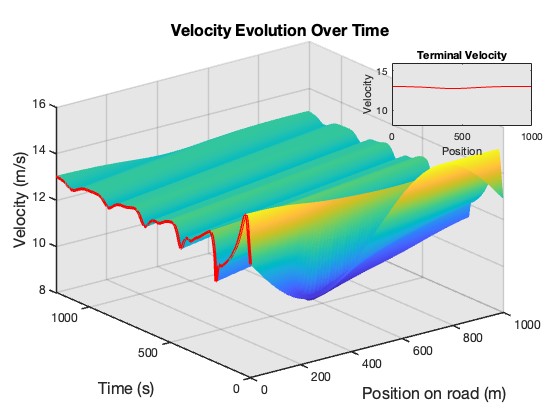}
        \caption{Velocity evolution}
        \label{212-2}
    \end{subfigure}
    \caption{Density and velocity evolution of the mixed flow model with $20\%$ of CAVs evenly distributed.}
    \label{Simu_even20}
\end{figure}

\begin{figure}[ht]
    \centering
        \centering
    \begin{subfigure}[b]{0.48\textwidth}
        \centering
        \includegraphics[height=0.2\textheight]{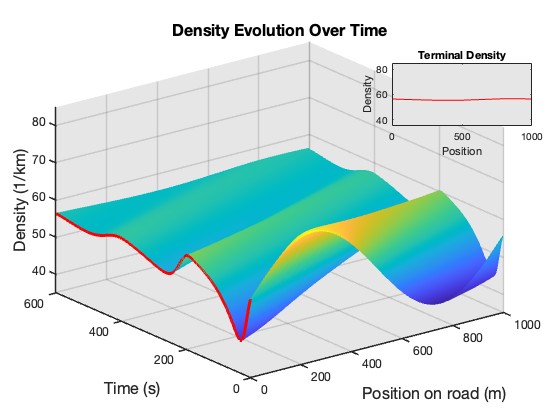}
        \caption{Density evolution}
        \label{213-1}
    \end{subfigure}
    \begin{subfigure}[b]{0.48\textwidth}
        \centering
        \includegraphics[height=0.2\textheight]{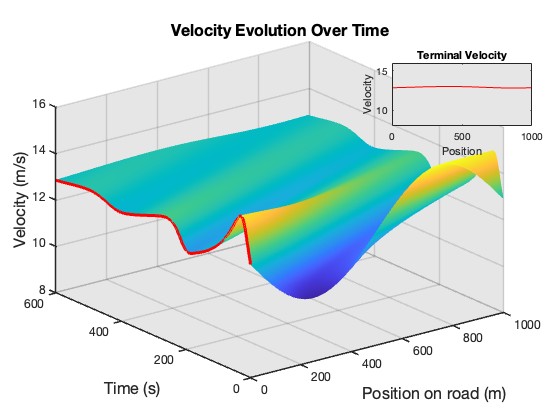}
        \caption{Velocity evolution}
        \label{213-2}
    \end{subfigure}
    \caption{Density and velocity evolution of the mixed flow model with $40\%$ of CAVs evenly distributed.}
    \label{Simu_even40}
\end{figure}

 We first consider CAVs evenly distributed in the mixed traffic with penetration rates of $10\%$, $20\%$, $40\%$. Based on the results of the single-class model, we choose the observation distance $L_D=100$m for CAVs. For $10\%$ and $20\%$ penetration rates we set the time duration to be $T=1200$s and for $40\%$ we set $T=600$s. The numerical results of density and velocity evolution are shown in \textbf{Figures \ref{Simu_even10}-\ref{Simu_even40}}.
 \begin{remark}
     For the mixed plot we plot the evolution of the total density and the HDVs' velocity, since traffic flow of pure CAVs are already shown stable.
 \end{remark}
% Blowup of shockwaves are more likely to happen when CAVs and HDVs are separated: in this case CAVs have less control of HDVs
% With proper front-following distance 20 percent of CAVs can stabilize the traffic
% Takes longer period of time to stabilize than previous models, but this is because the relaxation is not specifically designed for stabilization purpose, also the oscillation is bigger than Xuan Di's paper

%You have 10, 20 and 40$\%$ in the setting, not sure where this 5, 15, 45 numbers come from. 
%Those numbers are from different initial data. The long-term stability is not a good term and I have rewrite this part.
From these results, we can observe that $20$ percent of CAVs can stabilize the mixed flow to smaller oscillations, and $40$ percent of CAVs has faster convergence to equilibrium state, while $10$ percent of CAVs fails to stabilize the traffic. Such results are consistent with those from a similar study \cite{huang2020scalable}.

\subsubsection{Segregated distribution}
\begin{figure}[ht]
    \centering
        \centering
    \begin{subfigure}[b]{0.48\textwidth}
        \centering
        \includegraphics[height=0.2\textheight]{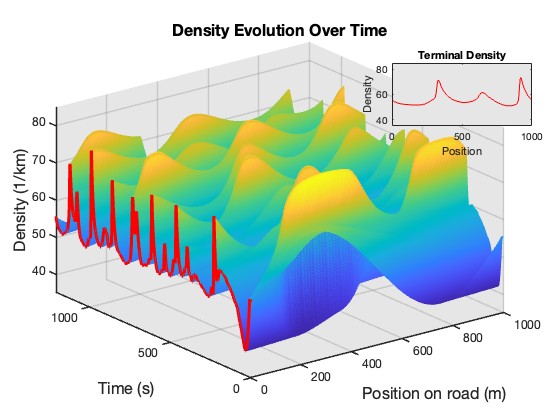}
        \caption{Density evolution}
        \label{221-1}
    \end{subfigure}
    \begin{subfigure}[b]{0.48\textwidth}
        \centering
        \includegraphics[height=0.2\textheight]{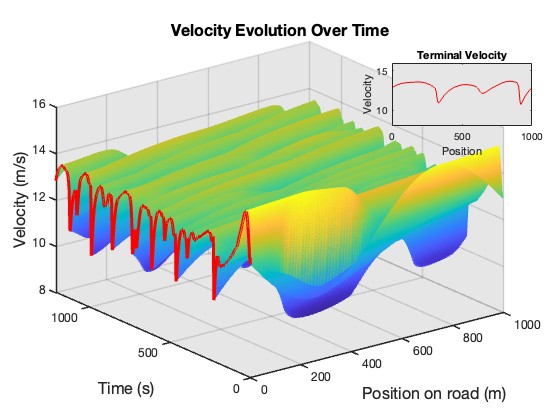}
        \caption{Velocity evolution}
        \label{221-2}
    \end{subfigure}
    \caption{Density and velocity evolution of the mixed flow model with $10\%$ of concentrated CAVs.}
    \label{Simu_con10}
\end{figure}

\begin{figure}[ht]
    \centering
        \centering
    \begin{subfigure}[b]{0.48\textwidth}
        \centering
        \includegraphics[height=0.2\textheight]{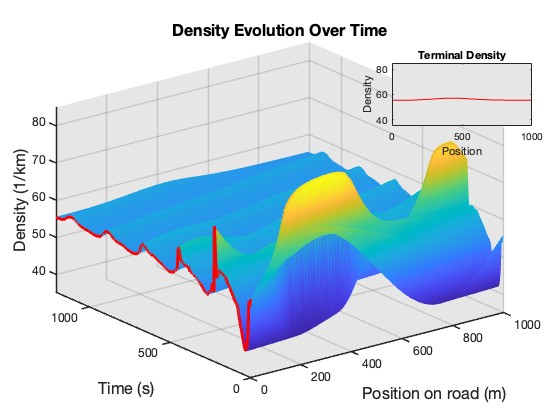}
        \caption{Density evolution}
        \label{222-1}
    \end{subfigure}
    \begin{subfigure}[b]{0.48\textwidth}
        \centering
        \includegraphics[height=0.2\textheight]{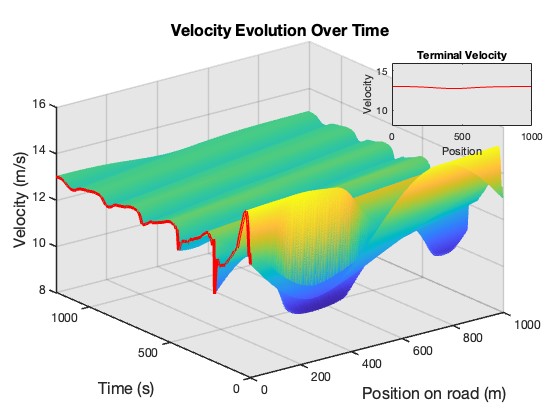}
        \caption{Velocity evolution}
        \label{222-2}
    \end{subfigure}
    \caption{Density and velocity evolution of the mixed flow model with $20\%$ of concentrated CAVs.}
    \label{Simu_con20}
\end{figure}

\begin{figure}[ht]
    \centering
        \centering
    \begin{subfigure}[b]{0.48\textwidth}
        \centering
        \includegraphics[height=0.2\textheight]{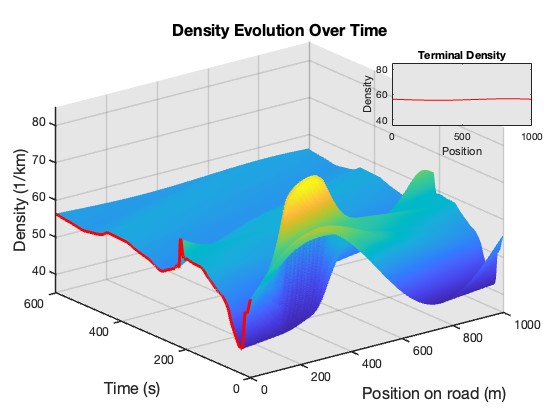}
        \caption{Density evolution}
        \label{223-1}
    \end{subfigure}
    \begin{subfigure}[b]{0.48\textwidth}
        \centering
        \includegraphics[height=0.2\textheight]{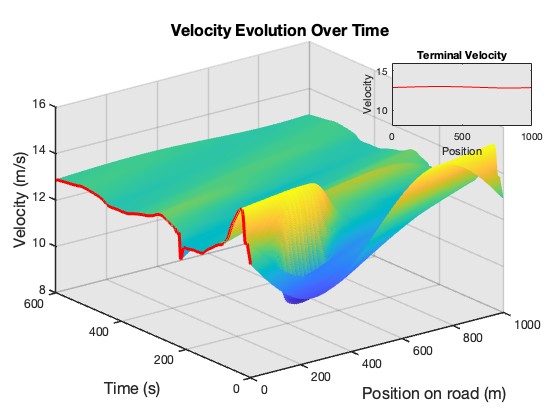}
        \caption{Velocity evolution}
        \label{223-2}
    \end{subfigure}
    \caption{Density and velocity evolution of the mixed flow model with $40\%$ of concentrated CAVs.}
    \label{Simu_con40}
\end{figure}

Now we consider another type of distribution such that CAVs and HDVs are segregated into two parts. With the same penetration rates, we let CAVs concentrate at around $x=500$m and HDVs concentrate at the remaining locations. In details, the initial density of CAVs is given as
\begin{equation}
    \rho^c(x)=\begin{cases}
        0.999\rho^s, \text{ if } \frac{1-r}{2}L < x < \frac{1+r}{2}L, \\
        0.001\rho^s, \text{ otherwise. }
    \end{cases},
\end{equation}
where $r$ is the percentage of CAVs, and $\rho^h$ can be calculated by $\rho^h=\rho^s-\rho^c$. The small densities is designed for numerical stability. We set the same simulation time as in the evenly distributed case. The numerical results of density and velocity evolution are shown in \textbf{Figures \ref{Simu_con10}-\ref{Simu_con40}}.


These results showed that the segregated distributions of mixed flow have similar asymptotic behaviors as even distributions. The main difference is that the initial waves have larger scales for segregated distributions where the HDVs are concentrated, since HDV traffic is less stable than CAV traffic. Overall the convergence of mixed traffic is slower than the results obtained in \cite{huang2020scalable}, possibly due to large oscillations and inadequate information utilized from HDVs. One possible improvement is to add predictive or feedback controls as previously investigated in car-following models e.g. \cite{zhou2019distributed, jin2020dynamical}.
\FloatBarrier
\section{Concluding remarks} \label{sec36}
%Consider a model as CAVs and gives the multiclass extension, analysed stability for different ability of CAVs and mixed traffic
%This is just an initial piece of work on second order multiclass traffic.
%Future work: 1 well-definedness for general model, extend model for CAVs, 2 consider more factor for HDVs and more vehicle classes, combine with other control methods, 3 combine with controls like ramp metering, and variable speed limit; 4 Evaluate performance of mixed traffic condition and road condition, extend to network model
% Contribution in model, stability, experiments
This paper make extensions to a second order non-equilibrium traffic flow model, i.e. the ARZ model, to take into account the look-ahead capabilities of CAVs, either in a single-class or multi-class context. The look-ahead effect is captured by a modification of the relaxation term, which can be interpreted as CAVs attempt to adopt a target speed based on the average traffic conditions within its spatial observation range, similar to multi-following microscopic traffic models. The stability properties of both extended models are analysed through wave perturbation analysis, and the results show that a longer observation range yields a less restrictive stability condition. Numerical solution using forward schemes with approximate Riemann solvers is provided,and numerical experiments are carried out to examine the effects of various parameters and the spatial distribution of CAVS on both the stability of mixed autonomy traffic, and CAVs' ability to stabilize mixed traffic flow. It is found that higher penetration rate of CAVs stabilize mixed traffic flow faster, which is consistent with similar studies in \cite{huang2020scalable}. 

% Superiority and suggestions
Our study reveals several new insights on mixed autonomy traffic. One interesting finding is that having more information of traffic conditions does not necessarily translate into better control of traffic. In our particular setting, a moderate look-ahead distance of $100$m enables faster convergence to equilibrium than having the full observation of road conditions on the entire ring road. Another interesting finding is that the distribution of vehicles have little effect on long-term stability of mixed traffic flow, but the initial oscillations for segregated distribution have larger amplitudes than that in the even distribution case. These insights can help CAV manufacturers design more effective control algorithm that can benefit both parties in mixed autonomy traffic, and traffic engineers to better manage mixed autonomy flow through leveraging the sensing and control capabilities of CAVs. 

% Future work
%This paper can give several promising future research directions, including: (1) Well-definedness and stability are not shown in details in this paper. More general boundary conditions and stability analysis still need to be learned (this is a very challenging topic). (2) CAVs can do more than just using average density of front-observation, as mentioned previously. Various other controls including delayed feedback control \cite{jin2020dynamical}, predictive control \cite{zhou2019distributed} can be implemented with CAVs and the interaction between CAVs and HDVs can be a more generalized form. (3) Other than using solely CAVs as controllers, we can combine with upper level controls, e.g. ramp metering \cite{yu2019traffic}, variable speed limit (VSL) \cite{yu2018varying} as cooperative controller to improve traffic condition. (4) Extensions with more general road conditions (multi-lane, lane-free) and more class of vehicles (trucks, motorcycles) are parts of future directions and can offer more insights for real-world traffic.

\bibliographystyle{unsrt}  
\bibliography{reference} 
\end{document}